\documentclass[twoside]{article} 
\usepackage{graphicx}
\date{} 
\title{The asymptotic expansion of the Humbert hyper-Bessel function}
\author{\sc R. B.\ Paris \\
{\em Division of Computing and Mathematics,} \\
{\em Abertay University, Dundee DD1 1HG, UK}}
\begin{document}
\def\f#1#2{\mbox{${\textstyle \frac{#1}{#2}}$}}
\def\dfrac#1#2{\displaystyle{\frac{#1}{#2}}}
\def\boldal{\mbox{\boldmath $\alpha$}}
\newcommand{\bee}{\begin{equation}}
\newcommand{\ee}{\end{equation}}
\newcommand{\lam}{\lambda}
\newcommand{\ka}{\kappa}
\newcommand{\al}{\alpha}
\newcommand{\ba}{\beta}
\newcommand{\la}{\lambda}
\newcommand{\ga}{\gamma}
\newcommand{\eps}{\epsilon}
\newcommand{\fr}{\frac{1}{2}}
\newcommand{\fs}{\f{1}{2}}
\newcommand{\g}{\Gamma}
\newcommand{\br}{\biggr}
\newcommand{\bl}{\biggl}
\newcommand{\ra}{\rightarrow}
\newcommand{\gtwid}{\raisebox{-.8ex}{\mbox{$\stackrel{\textstyle >}{\sim}$}}}
\newcommand{\ltwid}{\raisebox{-.8ex}{\mbox{$\stackrel{\textstyle <}{\sim}$}}}
\renewcommand{\topfraction}{0.9}
\renewcommand{\bottomfraction}{0.9}
\renewcommand{\textfraction}{0.05}
\newcommand{\mcol}{\multicolumn}
\date{}
\maketitle
\pagestyle{myheadings}
\markboth{\hfill \sc R. B.\ Paris  \hfill}
{\hfill \sc Humbert's function\hfill}
\begin{abstract}
We consider the asymptotic expansion of the Humbert hyper-Bessel function expressed in terms of a ${}_0F_2$ hypergeometric function by
\[J_{m,n}(x)=\frac{(x/3)^{m+n}}{m! n!}\,{}_0F_2(-\!\!\!-;m+1, n+1; -(x/3)^3)\]
as $x\to+\infty$, where $m$, $n$ are not necessarily non-negative integers. Particular attention is paid to the determination of the exponentially small contribution. The main approach utilised is that described by the author (J. Comput. Appl. Math. {\bf 234} (2010) 488-504); a leading-order estimate is also obtained by application of the saddle-point method applied to an integral representation containing a Bessel function. Numerical results are presented to demonstrate the accuracy of the resulting compound expansion.
\vspace{0.3cm}

\vspace{0.1cm}


\end{abstract}

\vspace{0.3cm}

\noindent $\,$\hrulefill $\,$

\vspace{0.3cm}

\begin{center}
{\bf 1.\ Introduction}
\end{center}
\setcounter{section}{1}
\setcounter{equation}{0}
\renewcommand{\theequation}{\arabic{section}.\arabic{equation}}
In \cite{PH}, P. Humbert introduced the function he termed a Bessel function of the third order (and which we call a hyper-Bessel function) by
\[J_{m,n}(x)=\frac{(x/3)^{m+n}}{m! n!}\,{}_0F_2(-\!\!\!-;m+1, n+1; -(x/3)^3)\]
\bee\label{e11}
\hspace{0.9cm}=\frac{(x/3)^{m+n}}{m! n!} \sum_{k=0}^\infty \frac{(-)^k (x/3)^{3k}}{(m+1)_k (n+1)_k k!}
\ee
in his development of operational calculus. Here ${}_0F_2$ is a hypergeometric function with two denominator parameters and $(a)_k=\g(a+k)/\g(a)$ is the Pochhammer symbol. He obtained several properties of this function and, in particular,
showed that
\[\sum_{n=0}^\infty \frac{(-\la x/3)^n}{n!}\,J_{n,n}(x)=J_{0,0}(x(1+\la)^{1/3});\]
when $\la=-1$, the right-hand side of this result reduces to $J_{0,0}(0)=1$. 

Further investigation into the convergence of
infinite series involving the above Humbert function was carried out by Varma \cite{V}. He considered series of the type
\[\sum_{n=1}^\infty a_n J_{m,n}(x)\qquad \mbox{and}\qquad \sum_{n=1}^\infty a_{n} J_{\al n+\beta,n}(x),\]
where $\al$ and $\beta$ are fixed positive constants.
The convergence domains of these series were obtained by determining the leading asymptotic behaviour of $J_{m,n}(x)$ for large $m$ and $n$. In a footnote in Varma's paper it was stated that the asymptotics of (\ref{e11}) for large $x$ would be determined, but a literature search has not revealed such a study.

An early paper by Wrinch \cite{DW} dealt with the asymptotics of the hypergeometric function ${}_0F_n(z)$ with $n$ denominator parameters. Using the assumed asymptotic form for the case of $n-1$ denominator parameters, she proceeded to employ an inductive argument to determine the asymptotic form in the case of $n$ parameters. Exponentially small contributions were discarded. Finally, we mention the paper \cite{PW} which considered the hyper-Bessel differential equation of order $n$
\[u^{(n)}-z^m u=0\]
and expressed solutions in terms of a variety of integral representations, both single and double.

The development of exponentially precise asymptotics during the past two decades has shown that retention of exponentially small expansions, although negligible in the Poincar\'e sense, can significantly improve the achievable numerical accuracy. A nice example that illustrates the advantage of retaining such terms is given in Olver's book \cite[p.~76]{Olv}.
In this paper we investigate the asymptotic expansion of the hyper-Bessel function $J_{m,n}(x)$ for $x\to+\infty$ using the well-established asymptotic theory of integral functions of hypergeometric type. We pay particular attention to the exponentially small contribution to (\ref{e11}), which we present in Section 3 following the approach employed in \cite{P}. An alternative approach using the saddle-point method to derive the leading terms applied to an integral representation of ${}_0F_2$ involving a Bessel function  is given in Section 4.

A numerical section displays results that confirm the accuracy of the expansion for the hyper-bessel function in (\ref{e11}). 
In Section 6, we present a summary of the asymptotic expansions of the function
\[\hspace{2cm}_0F_{n-1}(-\!\!\!-; b_1, \ldots , b_{n-1}; -(x/n)^n) \qquad (x\to+\infty)\]
for $n=4$ and $n=5$. In an appendix we give an algorithm for the determination of the coefficients appearing in the different expansions.

\vspace{0.6cm}

\begin{center}
{\bf 2.\ The dominant asymptotic expansion}
\end{center}
\setcounter{section}{2}
\setcounter{equation}{0}
\renewcommand{\theequation}{\arabic{section}.\arabic{equation}}
We replace $m$ and $n$ in (\ref{e11}) by $a-1$ and $b-1$ and consider the function
\bee\label{e21}
F(x):=\frac{1}{\g(a)\g(b)} \,{}_0F_2(-\!\!\!-;a,b;-(x/3)^3)=\sum_{k=0}^\infty \frac{(-)^k (x/3)^{3k}}{\g(k+a) \g(k+b) k!},
\ee
where $a$, $b$ are in general complex constants such that $a,\ b\neq 0, -1, -2, \ldots\ $. The function in (\ref{e21}) is a particular case of the general integral function of hypergeometric type given by
\bee\label{e22}
{}_p\Psi_q(z)=\sum_{n=0}^\infty g(n)\frac{z^n}{n!},\qquad g(n):= \frac{\prod_{r=1}^p \g(\al_r n+a_r)}{\prod_{r=1}^q \g(\beta_r n+b_r)},
\ee
where $p$ and $q$ are non-negative integers, the parameters $\alpha_r>0$,
$\beta_r>0$ and $a_r$ and $b_r$ are
arbitrary complex numbers. We also assume that the $\alpha_r$ and $a_r$ are subject to 
the restriction
\begin{equation}\label{e11b}
\alpha_rn+a_r\neq 0, -1, -2, \ldots \qquad (n=0, 1, 2, \ldots\ ;\, 1\leq r \leq p)
\end{equation}
so that no gamma function in the numerator in (\ref{e22}) is singular.

We introduce the parameters associated\footnote{Empty sums and products are to be interpreted as zero and unity, respectively.} with $g(n)$ in (\ref{e22}) given by
\[\kappa=1+\sum_{r=1}^q\beta_r-\sum_{r=1}^p\alpha_r, \qquad 
h=\prod_{r=1}^p\alpha_r^{\alpha_r}\prod_{r=1}^q\beta_r^{-\beta_r},\]
\begin{equation}\label{e12}
\vartheta=\sum_{r=1}^pa_r-\sum_{r=1}^qb_r+\f{1}{2}(q-p),\qquad \vartheta'=1-\vartheta.
\end{equation}
If it is supposed that $\alpha_r$ and $\beta_r$ are such that $\kappa>0$ then ${}_p\Psi_q(z)$ 
is uniformly and absolutely convergent for all finite $z$.
The parameter $\kappa$ plays a critical role 
in the asymptotic structure of ${}_p\Psi_q(z)$ by determining the sectors in the $z$-plane 
in which its behaviour is either exponentially large, algebraic or exponentially small 
in character as $|z|\ra\infty$. 

The exponential expansion ${E}(z)$ associated with ${}_p\Psi_q(z)$ is given by the {\it formal\/} asymptotic sum
\begin{equation}\label{e22c}
{E}(z):=Z^\vartheta e^Z\sum_{j=0}^\infty A_jZ^{-j}, \qquad Z=\kappa (hz)^{1/\kappa},
\end{equation}
where the coefficients $A_j$ are those appearing in the inverse factorial expansion of $g(s)/s!$ given by  
\begin{equation}\label{e22a}
\frac{g(s)}{\g(1+s)}=\kappa (h\kappa^\kappa)^s\bl\{\sum_{j=0}^{M-1}\frac{A_j}{\Gamma(\kappa s+\vartheta'+j)}
+\frac{\rho_M(s)}{\Gamma(\kappa s+\vartheta'+M)}\br\}.
\end{equation}
Here $g(s)$ is defined in (\ref{e22}) with $n$ replaced by $s$, $M$ is a positive integer and $\rho_M(s)=O(1)$ for $|s|\ra\infty$ in $|\arg\,s|<\pi$.
The leading coefficient $A_0$ is specified by
\begin{equation}\label{e22b}
A_0=(2\pi)^{\frac{1}{2}(p-q)}\kappa^{-\frac{1}{2}-\vartheta}\prod_{r=1}^p
\alpha_r^{a_r-\frac{1}{2}}\prod_{r=1}^q\beta_r^{\frac{1}{2}-b_r}.
\end{equation}
The coefficients $A_j$ are independent of $s$ and depend only on the parameters $p$, $q$, $\alpha_r$, 
$\beta_r$, $a_r$ and $b_r$. 
When $p\geq1$ there is also an algebraic expansion which we do not give here as the function in (\ref{e21}) corresponds to $p=0$ for which there is no algebraic expansion. 
The asymptotic expansion of ${}_p\Psi_q(z)$ when $\ka>2$ and $p=0$ is given by (see, for example, \cite[p,~58]{PK})
\bee\label{e23a}
{}_p\Psi_q(z)\sim \sum_{r=-P}^P E(ze^{2\pi ir})\qquad (|\arg\,z|\leq\pi),
\ee
where $P$ is chosen such that $2P+1$ is the smallest odd integer satisfying $2P+1>\fs\ka$.

For the function $F(x)$ in (\ref{e21}) we have $p=0$, $q=2$ and the parameters
\bee\label{e23b}
\ka=3,\quad h=1,\quad \vartheta=1-a-b,\quad \vartheta'=a+b,\quad A_0=\frac{3^{-\frac{1}{2}-\vartheta}}{2\pi},
\ee
and $P=1$, with $Z=xe^{\pi i/3}$ so that (\ref{e23a}) yields
\[F(x)\sim \sum_{r=-1}^1 E(xe^{(2r+1)\pi i/3})\qquad (x\to+\infty).\]
The dominant contribution results from the series with $r=0$ and $r=-1$ to yield
\[E(xe^{\pi i/3})+E(xe^{-\pi i/3})\hspace{10cm}\]
\[=(xe^{\pi i/3})^\vartheta e^{xe^{\pi i/3}} \sum_{j=0}^\infty A_j (xe^{\pi i/3})^{-j}+(xe^{-\pi i/3})^\vartheta e^{xe^{-\pi i/3}} \sum_{j=0}^\infty A_j (xe^{-\pi i/3})^{-j}\]
\bee\label{e23}
=2A_0 x^\vartheta e^{x/2} \sum_{j=0}^\infty c_j x^{-j} \cos \bl(\frac{\sqrt{3}}{2} x+\frac{\pi}{3} (\vartheta-j)\br).
\ee
The first few normalised coefficients $c_j:=A_j/A_0$ are given by (see the appendix) 
\[c_0=1,\qquad c_1=-\f{2}{3}+\wp_1-\wp_2+ab,\]
\bee\label{e2coeff}
c_2=\f{2}{9}+\f{1}{6}\{-4\wp_1+\wp_2-4\wp_3+3\wp_4-3ab(\wp_1+2\wp_2)+ab(17+9ab)\},
\ee
\[c_3=\f{32}{81}+\f{1}{162}\{-72\wp_1-198\wp_2+45\wp_3+81\wp_4+27\wp_5-27\wp_6+3ab(120+135ab+63a^2b^2)\]
\[+27a^2b^2(\wp_1-6\wp_2)+9ab(24\wp_1-63\wp_2+6\wp_3+9\wp_4)\},\]
where, for brevity, we have put $\wp_n:=a^n+b^n$. The rapidly increasing complexity of the higher coefficients prevents their explicit representation. However, when dealing with specific values of $a$ and $b$ it is possible to generate many coefficients; see Section 5. We observe that these coefficients are symmetrical in $a$ and $b$, which is necessary since the parameters $a$ and $b$ may be interchanged in (\ref{e21}). It can be verified that $c_1=c_2=c_3=0$ when $a=\f{1}{3}$, $b=\f{2}{3}$ and $a=\f{4}{3}$, $b=\f{5}{3}$; see (\ref{e26}), (\ref{e27}) and the appendix.

The expansion corresponding to $r=1$, namely 
\bee\label{e24}
E(xe^{\pi i})=A_0 e^{\pi i\vartheta} x^\vartheta e^{-x} \sum_{j=0}^\infty c_j (-x)^{-j}
\ee
is an exponentially small contribution. However, this cannot be the correct form, since for $x>0$ and real values of $a$ and $b$ such that $a+b$ is non-integer, (\ref{e24}) is complex-valued whereas $F(x)$ is real. It is tempting to add the expansion corresponding to $r=-2$ which yields the conjugate of (\ref{e24}) to produce the exponentially small contribution
\bee\label{e25}
2A_0 \cos \pi\vartheta \,x^\vartheta e^{-x} \sum_{j=0}^\infty c_j(-x)^{-j}.
\ee

The expression in (\ref{e25}) is also incorrect. This can be seen by inspecting the case $a=\f{1}{3}$, $b=\f{2}{3}$ ($\vartheta=0$) to find by the triplication formula for the gamma function
\bee\label{e2gamma}
\g(3z)=\frac{3^{3z-\frac{1}{2}}}{2\pi} \g(z) \g(z+\f{1}{3}) \g(z+\f{2}{3})
\ee
the evaluation
\begin{eqnarray}
F(x)&=&\sum_{k=0}^\infty \frac{(-)^k (x/3)^{3k}}{\g(k+\f{1}{3}) \g(k+\f{2}{3}) k!}=\frac{3^{1/2}}{2\pi} \sum_{k=0}^\infty \frac{(-x)^{3k}}{\g(3k+1)}\nonumber\\
&=&\frac{3^{-1/2}}{2\pi}\bl(2e^{x/2} \cos \frac{\sqrt{3}}{2}x +e^{-x}\br).\label{e26}
\end{eqnarray}
But (\ref{e23}) together with (\ref{e25}) in which $c_j=0$ ($j\geq 1$)  for $a=\f{1}{3}$, $b=\f{2}{3}$ (see the appendix) predict the result
\[\frac{3^{-1/2}}{2\pi}\bl(2e^{x/2} \cos \frac{\sqrt{3}}{2}x +2e^{-x}\br)\]
in which the subdominant term is {\it twice} the correct value.
Similarly, if we take $a=\f{4}{3}$, $b=\f{5}{3}$ ($\vartheta=-2$) with $c_j=0$ ($j\geq 1$) (see the appendix) we find
\bee\label{e27}
F(x)=\frac{3^{5/2}}{2\pi} \sum_{k=0}^\infty \frac{(-x)^{3k}}{\g(3k+3)}=\frac{3^{3/2}}{2\pi x^2}\bl(2e^{x/2} \cos \bl(\frac{\sqrt{3}}{2} x-\frac{2}{3}\pi\br)+e^{-x}\br)
\ee
and again the subdominant term in (\ref{e25}) is twice the correct value.

The treatment of the exponentially small contribution to the expansion of $F(x)$ is presented in Section 3. This is based on the approach employed in \cite{P}, which we repeat here to make the paper self-contained.

\vspace{0.6cm}

\begin{center}
{\bf 3.\ The exponentially small contribution to $F(x)$}
\end{center}
\setcounter{section}{3}
\setcounter{equation}{0}
\renewcommand{\theequation}{\arabic{section}.\arabic{equation}}
We write $F(x)$ in (\ref{e21}) as the Mellin-Barnes integral
\bee\label{e31}
F(x)=\frac{1}{2\pi i} \int_C \frac{\g(s) (x/3)^{-3s}}{\g(a-s) \g(b-s)}\,ds,
\ee
where $C$ denotes a loop described in the positive sense with endpoints at infinity in $\Re (s)<0$ that encloses all poles of $\g(s)$ at $s=0, -1, -2, \ldots\ $. Evaluation of the residues at these simple poles yields the sum in (\ref{e21}). From Stirling's formula $\g(z)\sim (2\pi)^{1/2} e^{-z} z^{z-1/2}$ for large $|z|$ in $|\arg\,z|<\pi$, the dominant behaviour of the modulus of the integrand as $|s|\to\infty$ is controlled by the factor $\exp [\ka \Re (s) \log\,|s|]$, so that the integral in (\ref{e31}) converges {\it without\/} restriction on $x$.

We rewrite the above integral as
\bee\label{e32}
F(x)=\frac{\pi^{-2}}{2\pi i}\int_C \g(s) \g(s+1-a) \g(s+1-b) G(s) (x/3)^{-3s} ds,
\ee
where
\[G(s)=\sin \pi(a-s) \sin \pi(b-s)=\fs\{\cos \pi(a-b)-\cos \pi(a+b-2s)\}.\]
Since there are no poles to the right of $C$, we may displace the path as far to the right as we please (but with endpoints at infinity still in $\Re (s)<0$), so that $|s|$ is everwhere large on the expanded loop. We can then employ the inverse factorial expansion \cite[Lemma 2.2, p.~39]{PK}
\bee\label{e33}
\g(s) \g(s+1-a) \g(s+1-b)=\frac{3^{1-3s}}{(2\pi)^{-2}}\bl\{\sum_{j=0}^{M-1} (-)^j A_j \g(3s+\vartheta-j)+\rho_M(3s) \g(3s+\vartheta-M)\br\}
\ee
valid as $|s|\to\infty$ in $|\arg\,s|<\pi$, where $M$ is a positive integer. The remainder function $\rho_M(3s)$
is analytic in $\Re (s)>0$ and satisfies $\rho_M(3s)=O(1)$ as $|s|\to\infty$ in $|\arg\,s|<\pi$. The coefficients $A_j$ are the same as those appearing in (\ref{e22a}) when $p=0$ (see \cite[p.~39]{PK}), with $A_0$ given in (\ref{e23b}). An algorithm for their determination is given in the appendix. 
Substitution of (\ref{e33}) into (\ref{e32})  (after replacement of $s$ by $s/3$) then leads to
\bee\label{e34}
F(x)=2\sum_{j=0}^{M-1} (-)^{j-1} A_j\,\frac{1}{2\pi i} \int_C x^{-s} \g(s+\vartheta-j) \{\cos \pi(a+b-\f{2}{3}s)-\cos \pi(a-b)\}\,ds+R_M,
\ee
where the remainder
\[R_M=\frac{2}{\pi i}\int_C \rho_M(s) \g(s+\vartheta-M) x^{-s} G(\f{1}{3}s)\,ds.\]
It is shown in \cite[Lemma 2.8, p.~72]{PK} that an order estimate for the above remainder integral is
O$(x^{\vartheta-M} e^{x/2})$ for the part of $G(s)$ containing $\cos \pi(a+b-2s)$ and O$(x^{\vartheta-M}e^{-x})$ for the part containing $\cos \pi(a-b)$.

We now make use of the Cahen-Mellin integral (see, for example, \cite[p.~90]{PK})
\bee\label{e35}
\frac{1}{2\pi i}\int_C \g(s+\alpha) z^{-s} ds=z^\alpha e^{-z}
\ee
valid for all $\arg\,z$ when $C$ is the same loop contour as in (\ref{e31}). Expressing $\cos \pi(a+b-\f{2}{3}s)$ in exponentials, we have for the factors $\exp [\mp\pi i(a+b-\f{2}{3}s)]$ the result
\[e^{\mp\pi i(a+b)} \sum_{j=0}^{M-1} (-)^{j-1} A_j\,\frac{1}{2\pi i}\int_C \g(s+\vartheta-j) (xe^{\mp 2\pi i/3})^{-s} ds
\]
\[=e^{\mp\pi i(a+b)} \sum_{j=0}^{M-1}(-)^{j-1} (xe^{\mp 2\pi i/3})^{\vartheta-j} \exp [-xe^{\mp 2\pi i/3}]=\exp [xe^{\pm\pi i/3}] \sum_{j=0}^{M-1} A_j (xe^{\pm\pi i/3})^{\vartheta-j}.\]
Thus we obtain the dominant contribution as $x\to+\infty$ given by
\bee\label{e36}
2x^\vartheta e^{x/2}\bl\{\sum_{j=0}^{M-1} A_j x^{-j} \cos \bl(\frac{\sqrt{3}}{2} x+\frac{\pi}{3}(\vartheta-j)\br)+O(x^{-M})\br\},
\ee
which is seen to correspond to (\ref{e23}) when we recall that the coefficients $c_j=A_j/A_0$.

The contribution to the integral in (\ref{e34}) from the factor $\cos \pi(a-b)$ is 
\[2\cos \pi(a-b) \sum_{j=0}^{M-1} (-)^j A_j\,\frac{1}{2\pi i} \int_C x^{-s} \g(s+\vartheta-j)\,ds=2\cos \pi(a-b)\,x^\vartheta e^{-x} \sum_{j=0}^{M-1} (-)^j A_jx^{-j}.\]
Then we have the exponentially small component of $F(x)$ given by
\bee\label{e37}
2A_0 \cos \pi(a-b) \, x^\vartheta e^{-x} \bl\{\sum_{j=0}^{M-1} (-)^j c_j x^{-j}+O(x^{-M})\br\}.
\ee

From (\ref{e36}) and (\ref{e37}) we have the following theorem (see also  \cite[Eq.~(16)]{YL}, \cite[p.~501]{P}):
\newtheorem{theorem}{Theorem}
\begin{theorem}$\!\!\!.$ \  Let $\vartheta=1-a-b$, $A_0=3^{-\frac{1}{2}-\vartheta}/(2\pi)$ and $M$ be a positive integer. Then the following expansion holds
\[\frac{1}{\g(a)\g(b)} \,{}_0F_2(-\!\!\!-;a,b;-(x/3)^3)\sim 2A_0x^\vartheta e^{x/2}\bl\{\sum_{j=0}^{M-1} c_j x^{-j} \cos \bl(\frac{\sqrt{3}}{2} x+\frac{\pi}{3}(\vartheta-j)\br)+O(x^{-M})\br\}\]
\[+2A_0 \cos \pi(a-b) \, x^\vartheta e^{-x} \bl\{\sum_{j=0}^{M-1} (-)^j c_j x^{-j}+O(x^{-M})\br\}\]
as $x\to+\infty$, where the first few coefficients $c_j$ are listed in (\ref{e2coeff}).
\end{theorem}
From Theorem 1 and (\ref{e11}) we then obtain the main result of the paper:
\begin{theorem}$\!\!\!.$ \ With $\vartheta=-m-n-1$, the Hunbert function $J_{m,n}(x)=(x/3)^{m+n} F(x)$
with $a=m+1$, $b=n+1$ $($where $m$, $n$ are not necessarily integers$)$. Then we have the expansion
\[J_{m,n}(x) \sim \frac{3^{1/2}}{\pi x}\,e^{x/2} \sum_{j=0}^\infty c_j x^{-j} \cos \bl(\frac{\sqrt{3}}{2} x+\frac{\pi}{3}(\vartheta-j)\br)\hspace{3cm}\]
\[\hspace{4cm}+\frac{3^{1/2}}{\pi x} \cos \pi(m-n)\,e^{-x} \sum_{j=0}^\infty (-)^jc_j x^{-j}\]
as $x\to+\infty$, where the first few coefficients $c_j$ are listed in (\ref{e2coeff}).
\end{theorem}

In Fig.~1 we show a plot of $e^{-x/2} J_{m.n}(x)$ for $m=\fs$ and $n=\f{2}{3}$. This is seen to possess an oscillatory structure with an algebraic decay (given by $x^{-1}$) superimposed. It is worth remarking that if the factorials $\g(k+m+1)$ and $\g(k+n+1)$ in (\ref{e21}) are absent, there is an extreme cancellation between the terms with the result that the sum equals $\exp\,[-(x/3)^3]$. The addition of a single gamma function in the denominator upsets this extreme cancellation to produce a Bessel function, which is oscillatory with an algebraic decay. The appearance of the second gamma function further destroys the cancellation between terms to result in $J_{m,n}(x)$ being an {\it exponentially growing} function controlled by $e^{x/2}$. Other perturbations of the negative exponential series to produce different behaviour have been considered in \cite{PExp}.

\begin{figure}[t]
	\begin{center} \includegraphics[width=0.6\textwidth]{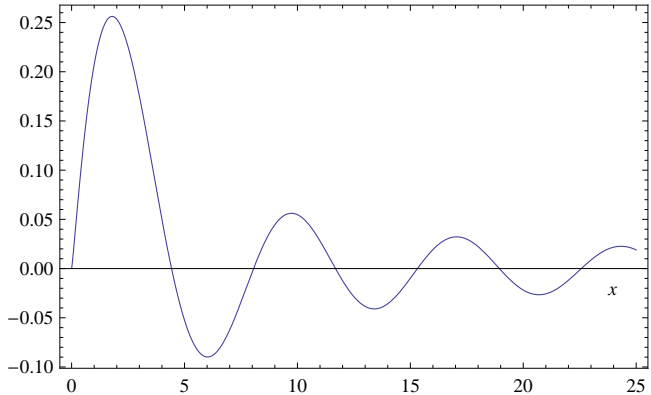}
\caption{\small{A plot of $e^{-x/2} J_{m,n}(x)$ when $m=1/2$ and $n=2/3$.}}
\end{center}
\end{figure}

\vspace{0.6cm}

\begin{center}
{\bf 4.\ An alternative approach}
\end{center}
\setcounter{section}{4}
\setcounter{equation}{0}
\renewcommand{\theequation}{\arabic{section}.\arabic{equation}}
In this section we give an alternative demonstration that the exponentially small contribution to $F(x)$ has the
coefficient $2A_0 \cos \pi(a-b)$ given in Theorem 1. We shall employ the saddle-point method to an integral representation for $F(x)$ and content ourselves with only the leading terms of the expansion.

We employ the Hankel loop integral for the gamma function \cite[(5.9.2)]{DLMF}
\[\frac{1}{\g(k+a)}=\frac{1}{2\pi i} \int_{-\infty}^{(0+)} \frac{e^t}{t^{a+k}}\,dt\]
in the series expansion for $F(x)$ in (\ref{e21}) to obtain
\[F(x)=\frac{1}{2\pi i}\int_{-\infty}^{(0+)} t^{-a}e^t \sum_{k=0}^\infty \frac{(-)^k ((x/3)^{3/2}/\sqrt{t})^{2k}}{\g(k+b) k!}\,dt\]
\bee\label{e41}
\hspace{3cm}=\frac{(x/3)^{3(1-b)/2}}{2\pi i} \int_{-\infty}^{(0+)} t^{\frac{1}{2}(b-1)-a}e^t\,J_{b-1}\bl(\frac{2(x/3)^{3/2}}{t^{1/2}}\br)\,dt,
\ee
where $J_\nu(z)$ denotes the Bessel function of the first kind. 

For $x>0$, the argument of the Bessel function appearing in (\ref{e41}) is in the range $[-\fs\pi,\fs\pi]$. As $|z|\to\infty$, we have the leading asymptotic behaviour \cite[(10.17.3)]{DLMF}
\[J_\nu(z)\sim \sqrt{\frac{2}{\pi z}}\,\cos (z-\fs\pi\nu-\f{1}{4}\pi)\qquad (|\arg\,z|<\pi).\]
Substitution of this asymptotic approximation into (\ref{e41}), followed by decomposition of the cosine into exponentials and the change of variable $t\to x\tau/3$,  leads to the approximation
\[F(x)\sim \frac{(x/3)^{\vartheta+\frac{1}{2}}}{2\sqrt{\pi}}\,\bl\{e^{-\frac{1}{2}\pi i(b-\frac{1}{2})}I_++e^{\frac{1}{2}\pi i(b-\frac{1}{2})}I_-\br\},\]
where
\bee\label{e42}
I_\pm=
\,\frac{1}{2\pi i} \int_{-\infty}^{(0+)} \tau^\gamma \exp\,[(x/3) \psi_\pm(\tau)]\, d\tau,\qquad \psi_\pm(\tau)=\tau\pm \frac{2i}{\sqrt{\tau}}
\ee
and $\gamma=\fs b-a-\f{1}{4}$. We shall find that the principal contributions to $I_\pm$ as $x\to+\infty$ arise from a neighbourhood of the unit circle in the $\tau$-plane (that is, of the circle of radius $x$ in the $t$-plane), thereby ensuring the validity of the use of the asymptotic approximation for the Bessel function in (\ref{e41}) since its argument is $O(x)$.

\begin{figure}[th]
	\begin{center} \includegraphics[width=0.4\textwidth]{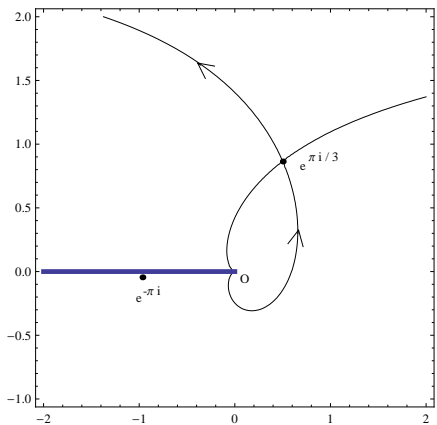}
\caption{\small{The paths of steepest descent through the saddle points $e^{\pi i/3}$ and $e^{-\pi i}$ for the phase function $\psi_+(\tau)$. The path through $e^{-\pi i}$ coincides with the lower side of the branch cut on the negative real $\tau$-axis. The arrows indicate the direction of integration over the saddle at $e^{\pi i/3}$.}}
\end{center}
\end{figure}

We consider the integral $I_+$. Saddle points of $\psi_+(\tau)$ occur at $\tau^{3/2}=i$; that is, at the points 
$\tau_{s1}=e^{\pi i/3}$ and $\tau_{s2}=e^{-\pi i}$. The path of steepest descent through $e^{\pi i/3}$ emanates from the origin as shown in Fig.~2 and passes to infinity in the second quadrant; the steepest descent path through $e^{-\pi i}$ lies on the lower side of the branch cut situated on the negative real axis. The loop integral can be made to coincide with these two paths. Then a straightforward application of the saddle-point method applied to the saddle $\tau_{s1}$, with $\psi_+(\tau_{s1})=3e^{\pi i/3}$ and $\psi_+''(\tau_{s1})=(3/2) e^{-\pi i/3}$ (so that the direction of integration through $\tau_{s1}$ is $2\pi/3$), yields the contribution
\[\frac{e^{\frac{1}{2}\pi i(b-\frac{1}{2})}}{2\pi i}\, \cdot\, 2\sqrt{\frac{\pi}{x}}\,\exp\,[xe^{\pi i/3}]\,e^{\pi i\gamma/3}\,e^{\pi i/3}=\frac{\exp\,[xe^{\pi i/3}]}{\sqrt{\pi x}}\,e^{\pi i\vartheta/3}.\]
A similar treatment for the integral $I_-$, where the steepest descent paths in Fig.~2 are replaced by the conjugate
paths with the saddles now situated at $e^{-\pi i/3}$ and $e^{\pi i}$ on the upper side of the branch cut. This yields the conjugate expression, so that the leading form of the dominant contribution to $F(x)$ is
\[2A_0 x^\vartheta e^{x/2}\,\cos \bl(\frac{\sqrt{3}}{2} x+\frac{\pi\vartheta}{3}\br)\qquad (x\to+\infty),\]
which agrees with that given in Theorem 1.

A similar treatment for the saddle $\tau_{s2}=e^{-\pi i}$, where $\psi_+(\tau_{s2})=-3$ and $\psi_+''(\tau_{s2})=-3/2$, yields the contribution to $e^{-\frac{1}{2}\pi i(b-\frac{1}{2})}I_+$ given by
\[\frac{e^{-\frac{1}{2}\pi i(b-\frac{1}{2})}}{2\pi i}\,\cdot\,2\sqrt{\frac{\pi}{x}}\,e^{-x-\pi i\gamma}=
\frac{e^{-x}}{\sqrt{\pi x}}\,e^{\pi i(a-b)}.\]
With the conjugate expression arising from the integral $e^{\frac{1}{2}\pi i(b-\frac{1}{2})}I_-$, we therefore obtain the leading-order subdominant contribution to $F(x)$ given by
\[2A_0 \cos \pi(a-b)\,x^\vartheta e^{-x}\qquad (x\to+\infty)\]
in accordance with the result stated in Theorem 1.
\vspace{0.6cm}

\begin{center}
{\bf 5.\ Numerical results}
\end{center}
\setcounter{section}{5}
\setcounter{equation}{0}
\renewcommand{\theequation}{\arabic{section}.\arabic{equation}}
To compute the series in Theorem 1 to high accuracy necessitates the evaluation of the normalised coefficients $c_j$ ($j\geq 1$). An algorithm for the calculation of these coefficients is summarised in the appendix, where in our calculations we have employed up to 25 coefficients. We present in Table 1 the coefficients with $1\leq j\leq 10$ for the case of $F(x)$ with $a=2/3$ and $b=5/6$.
\begin{table}[t]
\caption{\footnotesize{The normalised coefficients $c_j$ for $1\leq j\leq 10$ (with $c_0=1$) for $a=2/3$ and $b=5/6$.}}
\begin{center}
\begin{tabular}{|cc|cc|}
\hline
\mcol{1}{|c}{$j$} & \mcol{1}{c|}{$c_j$} & \mcol{1}{c}{$j$}& \mcol{1}{c|}{$c_j$} \\
[.1cm]\hline
&&&\\[-0.25cm]
1 & $+0.250000000000000$ & 6 & $-2.249984741210938$ \\
2 & $+0.156250000000000$ & 7 & $-9.259891510009766$ \\
3 & $+0.117187500000000$ & 8 & $-33.11664164066314$ \\
4 & $+0.017089843750000$ & 9 & $-90.94397798180580$ \\
5 & $-0.422973632812500$ & 10 & $-18.51875754073262$ \\
[.1cm]\hline
\end{tabular}
\end{center}
\end{table}

In Table 2 we show values of the absolute relative error in the computation of $F(x)$ in (\ref{e21}) using 
the optimally truncated\footnote{Optimal truncation corresponds to truncation of the asymptotic series at, or near, the term of least magnitude.} exponentially large expansion in Theorem 1 for different values of $x$.
\begin{table}[t]
\caption{\footnotesize{The absolute relative error in the computation of $F(x)$ using the optimally truncated dominant expansion in Theorem 1 for different $x$ and parameters $a$ and $b$.}}
\begin{center}
\begin{tabular}{|c|c|c|c|}
\hline
&&&\\[-0.3cm]
\mcol{1}{|c|}{} & \mcol{1}{c|}{$a=2/3,\ b=5/6$} &\mcol{1}{c|}{$\ \ \ a=1,\ b=1\ \ \ $}& \mcol{1}{c|}{$\ \ a=3/2,\ b=1\ \ $} \\
\mcol{1}{|c|}{$x$} & \mcol{1}{c|}{Relative Error} & \mcol{1}{c|}{Relative Error}& \mcol{1}{c|}{Relative Error} \\
[.1cm]\hline
&&&\\[-0.25cm]
10 & $9.801\times 10^{-07}$ & $2.328\times 10^{-06}$ & $4.158\times 10^{-09}$ \\
15 & $1.363\times 10^{-10}$ & $1.048\times 10^{-08}$ & $7.850\times 10^{-11}$ \\
20 & $1.655\times 10^{-13}$ & $5.425\times 10^{-11}$ & $8.025\times 10^{-13}$ \\
25 & $5.649\times 10^{-17}$ & $3.233\times 10^{-13}$ & $2.214\times 10^{-15}$ \\
30 & $1.694\times 10^{-19}$ & $1.135\times 10^{-15}$ & $2.895\times 10^{-16}$ \\
[.1cm]\hline
\end{tabular}
\end{center}
\end{table}
To detect the presence of the exponentially small expansion present in $F(x)$, we compare the values of
\[{\cal F}(x):=F(x)- 2A_0x^\vartheta e^{x/2}\sum_{j=0}^{j_0} c_j x^{-j} \cos \bl(\frac{\sqrt{3}}{2} x+\frac{\pi}{3}(\vartheta-j)\br),\]
where $j_0$ denotes the optimal truncation index, with the values of the exponentially small expansion
\[{\cal E}_s(x):=2A_0 \cos \pi(a-b) \, x^\vartheta e^{-x} \sum_{j=0}^{\infty} (-)^j c_j x^{-j}.\]
The results in Table 3 appear to confirm the form of the exponentially small contribution to $F(x)$ in Theorem 1.
We remark that when $a-b$ equals half-integer values, ${\cal E}_s(x)$ vanishes.
\begin{table}[h]
\caption{\footnotesize{Values of ${\cal F}(x)$ and ${\cal E}_s(x)$ for different $x$ and two values of the parameters $a$ and $b$.}}
\begin{center}
\begin{tabular}{|c|c|c|c|}
\hline
&&&\\[-0.3cm]
\mcol{1}{|c|}{$a=2/3,\ b=4/3$} & \mcol{1}{|c|}{$x=10,\ j_0=12$} & \mcol{1}{c|}{$x=15,\ j_0=17$} & \mcol{1}{c|}{$x=20,\ j_0=24$} \\
[.1cm]\hline
&&&\\[-0.25cm]
${\cal F}(x)$ & $-4.43157\times 10^{-06}$ & $-2.23785\times 10^{-08}$ & $-1.24121\times 10^{-10}$ \\
${\cal E}_s(x)$ & $-4.21754\times 10^{-06}$ & $-2.20456\times 10^{-08}$ & $-1.24469\times 10^{-10}$ \\
[.1cm]\hline
&&&\\[-0.3cm]
\mcol{1}{|c|}{$a=5/4,\ b=1/4$} & \mcol{1}{|c|}{$x=10,\ j_0=13$} & \mcol{1}{c|}{$x=15,\ j_0=15$} & \mcol{1}{c|}{$x=20,\ j_0=24$} \\
[.1cm]\hline
&&&\\[-0.25cm]
${\cal F}(x)$ & $-4.87509\times 10^{-06}$ & $-2.57286\times 10^{-08}$ & $-1.50868\times 10^{-10}$ \\
${\cal E}_s(x)$ & $-4.77620\times 10^{-06}$ & $-2.59124\times 10^{-08}$ & $-1.50113\times 10^{-10}$ \\
[.1cm]\hline
\end{tabular}
\end{center}
\end{table}

\vspace{0.6cm}

\begin{center}
{\bf 6.\ Extension of $F(x)$}
\end{center}
\setcounter{section}{6}
\setcounter{equation}{0}
\renewcommand{\theequation}{\arabic{section}.\arabic{equation}}
In this section we present the expansions for $n=4$ and $n=5$ of the extended function
\[F_n(x):=\frac{1}{\prod_{j=1}^{n-1}\g(b_j)}\,{}_0F_{n-1}(-\!\!\!-; b_1, \ldots , b_{n-1}; -(x/n)^n)\]
\bee\label{e61}
=\sum_{k=0}^\infty \frac{(-)^k (x/n)^{nk}}{\g(k+b_1)\ldots \g(k+b_{n-1}) k!}
\ee
consisting of $n-1$ denominator parameters $b_1, \ldots , b_{n-1}$; when $n=3$, this function is equivalent to that in (\ref{e21}). The parameters associated with $F_n(x)$ are
\[\kappa=n,\qquad \vartheta=\frac{n-1}{2}-\sigma_n,\qquad \sigma_n=\sum_{j=1}^{n-1}b_j,\qquad A_0=\frac{n^{-\frac{1}{2}-\vartheta}}{(2\pi)^{(n-1)/2}}.\]

From (\ref{e31}) and (\ref{e32}) we have
\begin{eqnarray}
F_n(x)&=&\frac{1}{2\pi i} \int_C \frac{\g(s) (x/n)^{-ns}}{\prod_{j=1}^{n-1} \g(b_j-s)}\,ds\nonumber\\
&=& \frac{\pi^{1-n}}{2\pi i} \int_C \g(s) \prod_{j=1}^{n-1} \g(s\!+\!1\!-\!b_j)\,G_n(s)\,(x/n)^{-ns} ds,\label{e62}
\end{eqnarray}
where
\[G_n(s):=\prod_{j=1}^{n-1} \sin \pi(b_j-s).\]
As in Section 3, the contour $C$ can be displaced as far to the right as we please (but with endpoints still in $\Re (s)<0$) so that we can invoke the inverse factorial expansion \cite[Lemma 2.2, p.~39]{PK}
\bee\label{e63}
\g(s) \prod_{j=1}^{n-1} \g(s\!+\!1\!-\!b_j)\sim \frac{n^{1-ns}A_0}{(2\pi)^{1-n}} \sum_{j=0}^\infty (-)^j c_j \g(ns+\vartheta-j)\qquad (|s|\to\infty,\ |\arg\,s|<\pi)
\ee
where \cite[Eq.~(A.7)]{P}
\[c_0=1,\qquad c_1=\frac{n}{2}\bl\{\sum_{j=1}^{n-1} b_j(1-b_j)-\frac{\vartheta(1-\vartheta)}{n}-\frac{(n^2-1)}{6n}\br\}.\]
The coefficients $c_j$ ($j\geq 2$) can be obtained by the algorithm outlined in the appendix.
\vspace{0.3cm}

\noindent{\bf The case $n=4$.}
\vspace{0.2cm}

\noindent We have after some routine algebra
\[G_4(s)=\prod_{j=1}^3 \sin \pi(b_j-s)=\frac{1}{4}\bl\{\cos \pi(\vartheta+3s)-\sum_{j=1}^3 \cos \pi(\vartheta+2b_j+s)\br\}.\]
Substitution of this expression into the integral in (\ref{e62}) together with the expansion (\ref{e63}) then yields
\[
F_4(x) \sim 2A_0 x^\vartheta e^{x/\sqrt{2}}\, \sum_{j=0}^\infty c_j x^{-j} \cos \bl(\frac{x}{\sqrt{2}} +\frac{\pi}{4}(\vartheta-j)\br)\hspace{4cm}\]
\bee\label{e64}
\hspace{1cm}-2A_0 x^\vartheta e^{-x/\sqrt{2}} \sum_{j=0}^\infty c_j x^{-j} \sum_{r=1}^3 \cos \bl(\frac{x}{\sqrt{2}}+\frac{3\pi}{4}(\vartheta-j)+2\pi b_r\br)
\ee
as $x\to+\infty$, where $A_0=4^{-\frac{1}{2}-\vartheta}/(2\pi)^{3/2}$, $\vartheta=\f{3}{2}-\sigma_4$. The exponentially small expansion can be written alternatively as
\bee\label{e64a}
-2A_0 x^\vartheta e^{-x/\sqrt{2}} \sum_{j=0}^\infty c_j x^{-j} \bl\{\cos \phi \sum_{r=1}^3 \cos 2\pi b_r- \sin \phi \sum_{r=1}^3 \sin 2\pi b_r\br\},
\ee
where $\phi:=x/\sqrt{2}+3\pi(\vartheta-j)/4$,
and is seen to agree with that given\footnote{There is an obvious misprint in \cite[Eq.~(18)]{YL}.} in \cite[Eq.~(18)]{YL}. 

In the special case $b_j=j/4$ ($1\leq j\leq 3$), the gamma functions in the denominator of (\ref{e61}) combine to yield $(2\pi)^{3/2}4^{-\frac{1}{2}-4k} \g(4k+1)$ to produce the exact evaluation
\[F_4(x)=\frac{1}{2^{1/2} \pi^{3/2}} \sum_{k=0}^\infty \frac{(-)^k x^{4k}}{\g(4k+1)}=4A_0 \cos \frac{x}{\sqrt{2}}\,\cosh \frac{x}{\sqrt{2}}~.\] 
Since $\sum_{r=1}^3 \cos 2\pi b_r=-1$, $\sum_{r=1}^3 \sin 2\pi b_r=0$, $\vartheta=0$ and $c_j=0$ ($j\geq 1$), the right-hand side of the expansion (\ref{e64}) yields
\[2A_0 e^{x/\sqrt{2}} \cos \frac{x}{\sqrt{2}}+2A_0 e^{-x/\sqrt{2}} \cos \frac{x}{\sqrt{2}}=4A_0 \cos \frac{x}{\sqrt{2}}\,\cosh \frac{x}{\sqrt{2}}.\] 

\vspace{0.3cm}

\noindent{\bf The case $n=5$.}
\vspace{0.2cm}

\noindent We have (see \cite[p.~50]{Hob})
\[G_5(s)=\prod_{j=1}^4 \sin \pi(b_j-s)=\frac{1}{8}\bl\{\cos \pi(\vartheta+4s)-\sum_{j=1}^4 \cos \pi(\vartheta+2b_j+2s)+\sum_{j=1}^3 \cos \pi(\vartheta+2b_j+2b_4)\br\}\]
with $\vartheta=2-\sum_{j=1}^4 b_j$. Then the same procedure yields the expansion
\[
F_5(x)\sim 2A_0 x^\vartheta e^{x \cos \pi/5} \sum_{j=0}^\infty c_j x^{-j} \cos \bl(x\sin \frac{\pi}{5}+\frac{\pi}{5} (\vartheta-j)\br)\]
\[\hspace{3.2cm}-2A_0x^\vartheta e^{x \cos 3\pi/5} \sum_{j=0}^\infty c_j x^{-j}\bl\{\cos \phi \sum_{r=1}^4 \cos 2\pi b_r-\sin \phi \sum_{r=1}^4\sin 2\pi b_r\br\}\]
\bee\label{e65}
\hspace{0.8cm}+2A_0 \sum_{r=1}^3 \cos \pi (\vartheta+2b_r+2b_4) \,x^\vartheta e^{-x} \sum_{j=0}^\infty (-)^j c_j x^{-j}
\ee
as $x\to+\infty$, where $\phi:=x\sin 3\pi/5+3\pi(\vartheta-j)/5$, $A_0=5^{-\frac{1}{2}-\vartheta}/(2\pi)^2$ and $\vartheta=2-\sigma_5$.

In the special case $b_j=j/5$ ($1\leq j\leq 4$), we have
\begin{eqnarray}
F_5(x)&=&\frac{5^{-1/2}}{(2\pi)^2} \sum_{k=0}^\infty \frac{(-)^k x^{5k}}{\g(5k+1)}\nonumber\\
&=&\frac{5^{-1/2}}{(2\pi)^2} \bl\{2e^{x \cos \pi/5} \cos (x \sin \pi/5)+2e^{x \cos 3\pi/5} \cos (x\sin 3\pi/5)+e^{-x}\br\}.\label{e66}
\end{eqnarray}
Since $\vartheta=0$, $c_j=0$ ($j\geq 1$) and
\[\sum_{r=1}^4 \cos 2\pi b_r=-1,\qquad \sum_{r=1}^4 \sin 2\pi b_r=0, \qquad \sum_{r=1}^3 \cos \pi (\vartheta+2b_r+2b_4)=\fs,\]
it is seen that the right-hand side of (\ref{e65}) correctly reduces to the exact expression (\ref{e66}).

In Table 4 we show the values of
\[{\cal F}_4(x):=F_4(x)-2A_0 x^\vartheta e^{x/\sqrt{2}}\, \sum_{j=0}^\infty c_j x^{-j} \cos \bl(\frac{x}{\sqrt{2}} +\frac{\pi}{4}(\vartheta-j)\br),\]
where the dominant expansion is optimally truncated at index $j_0$,
compared with the exponentially small expansion in (\ref{e64a})
\[{\cal E}_s(x):=-2A_0 x^\vartheta e^{-x/\sqrt{2}} \sum_{j=0}^\infty c_j x^{-j} \bl\{\cos \phi \sum_{r=1}^3 \cos 2\pi b_r- \sin \phi \sum_{r=1}^3 \sin 2\pi b_r\br\},\]
where $\phi:=x/\sqrt{2}+3\pi(\vartheta-j)/4$. It is seen that there is reasonable agreement between ${\cal F}_4(x)$ and ${\cal E}_s(x)$ thereby lending support
to the result stated in (\ref{e64}).
\begin{table}[h]
\caption{\footnotesize{Values of ${\cal F}_4(x)$ and ${\cal E}_s(x)$ for different $x$ and two sets of values of the parameters.}}
\begin{center}
\begin{tabular}{|c|c|c|}
\hline
&&\\[-0.3cm]
\mcol{1}{c}{} & \mcol{1}{|c|}{$a=-\f{1}{4},\ b=\fs,\ c=\f{5}{8}$} & \mcol{1}{|c|}{$a=\f{3}{4},\ b=\f{4}{5},\ c=\fs$} \\
\mcol{1}{|c|}{} & \mcol{1}{|c|}{$x=15,\ j_0=17$} & \mcol{1}{c|}{$x=18,\ j_0=17$} \\
[.1cm]\hline
&&\\[-0.25cm]
${\cal F}_4(x)$ & $+8.51145\times 10^{-06}$ & $+1.65559\times 10^{-07}$  \\
${\cal E}_s(x)$ & $+8.56645\times 10^{-06}$ & $+1.67468\times 10^{-07}$  \\
[.1cm]\hline
\end{tabular}
\end{center}
\end{table}
\vspace{0.6cm}

\begin{center}
{\bf Appendix A: \ The coefficients $c_j$}
\end{center}
\setcounter{section}{1}
\setcounter{equation}{0}
\renewcommand{\theequation}{\Alph{section}.\arabic{equation}}
The coefficients $A_j$ appearing in the inverse factorial expansion (\ref{e33}) are the same as those in (\ref{e22a}) when $p=0$, $q=2$ and the parameters $\beta_1=\beta_2=1$, $b_1=a$, $b_2=b$ in (\ref{e22}); see \cite[p.~39]{PK}. Thus, it is sufficient to consider the inverse factorial expansion
\bee\label{a0}
\frac{1}{\g(s+a)\g(s+b)\g(s+1)}=3^{3s+1} \bl\{\sum_{j=0}^{M-1} \frac{A_j}{\g(3s+\vartheta'+j)}+\frac{O(1)}{\g(3s+\vartheta'+M)}\br\}
\ee
for positive integer $M$, where we recall that $\vartheta'=a+b$ and $A_0$ is given in (\ref{e23b}). This expansion may be rewritten as
\bee\label{a1}
\frac{\g(3s+\vartheta')}{\g(s+a)\g(s+b)\g(s+1)}=3^{3s+1}A_0 \bl\{\sum_{j=0}^{M-1} \frac{c_j}{(3s+\vartheta')_j}+\frac{O(1)}{(3s+\vartheta')_M}\br\},
\ee
where $c_j=A_j/A_0$. The algorithm based on (\ref{a1}) relies on the asymptotic expansion of the gamma function to determine recursively the coefficients $c_j$. This has been described in \cite[Appendix]{P}, \cite[\S 2.2.4]{PK}, 
and will not be repeated here.

In the special case $a=\f{1}{3}$, $b=\f{2}{3}$ ($\vartheta'=1$) we have from (\ref{e2gamma}) that the left-hand side of (\ref{a0}) equals $3^{3s+\frac{1}{2}}/(2\pi \g(3s+1))$, whence it follows that $c_0=1$ and $c_j=0$ ($j\geq 1$). Similarly, when $a=\f{4}{3}$, $b=\f{5}{3}$ ($\vartheta'=3$), the left-hand side of (\ref{a0}) becomes 
$3^{3s+\frac{5}{2}}/(2\pi \g(3s+3))$ and again $c_0=1$, $c_j=0$ ($j\geq 1$).

An alternative method of determining the coefficients $c_j$ in the form of a recurrence relation is given in \cite[\S 2.2.2]{PK} based on the paper by Riney \cite{R}. This takes the form
\bee\label{a2}
c_j=-\frac{1}{27j} \sum_{k=0}^{j-1} c_k \,e(j,k),\qquad e(j,k)=\sum_{r=1}^3 D_r \,\frac{(\vartheta'-3b_r)_{3+j}}{(\vartheta'-3b_r)_k}
\ee
where $b_1=a$, $b_2=b$ and $b_3=1$ and the coefficients $D_r$ are given by
\[D_1=-\frac{1}{(a-b)(1-a)},\quad D_2=\frac{1}{(a-b)(1-b)},\quad D_3=\frac{1}{(1-a)(1-b)}.\]

A disadvantage of (\ref{a2}) is the fact that the coefficients $D_r$ present singularities when $a=b$ and $a$ (and/or $b$) $=1$. This then requires a limiting procedure. The resulting coefficients $c_j$, however, are non-singular; see the first few given in (\ref{e2coeff}). This problem does not arise with the algorithm based on (\ref{a1}).

\vspace{0.6cm}

\noindent{ {\bf Acknowledgement:} \, The author wishes to acknowledge F. Mainardi for bringing to his attention the paper by Humbert and for suggesting an investigation into the asymptotic behaviour of $J_{m,n}(x)$.  


\vspace{0.6cm}

\begin{thebibliography}{99}
\bibitem{Hob}
E.W. Hobson, {\it A Treatise on Plane Trigonometry}, Cambridge University Press, Cambridge, 1925.

\bibitem{PH}
P. Humbert, Nouvelles remarques sur les fonctions de Bessel du troisi{\`e}me ordre, Atti. Pont. Accad. della Scienza, {\bf 87} (1933) 323--331.

\bibitem{YL}
Y.L. Luke, {\it Special Functions and Their Approximations}, Vol.~1, Academic Press, New York, 1969.

\bibitem{Olv}
F.W.J. Olver, {\it Asymptotics and Special Functions}, Academic Press, New York, 1974; Reprinted in A.K. Peters, Massachussets, 1997.

\bibitem{DLMF}
F.W.J. Olver, D.W. Lozier, R.F. Boisvert and C.W. Clark (eds.),    
{\it NIST Handbook of Mathematical Functions}, Cambridge University Press, Cambridge, 2010.

\bibitem{PExp}
R.B. Paris, On the growth of a class of perturbation of the exponential series. 
Math. Balkanica, {\bf 21} (2007) 183--200.

\bibitem{P} 
R.B. Paris, Exponentially small expansions in the asymptotics of the Wright function, J. Comput. Appl. Math., {\bf 234} (2010) 488--504.

\bibitem{PK}
R.B. Paris and D. Kaminski, {\em Asymptotics and Mellin-Barnes Integrals}, Cambridge University Press, Cambridge, 2001.

\bibitem{PW}
R.B. Paris and A.D. Wood, Results old and new on the hyper-Bessel equation,  Proc. Roy. Soc. Edinburgh, {\bf 106A} (1987) 259--265.

\bibitem{R} 
T.D. Riney, On the coefficients in asymptotic factorial expansions, Proc. Amer. Math. Soc., {\bf 7} (1956) 245--249.

\bibitem{V}
R.S. Varma, On Humbert functions, Ann. Math., {\bf 42} (1941) 429--436.

\bibitem{DW}
D.M. Wrinch, A generalized hypergeometric function with $n$ parameters, Phil. Mag., {\bf 41} (1921) 174--186.
\end{thebibliography}
\end{document}